\documentclass[12pt]{amsart}
\usepackage{amsmath}
\usepackage{amsthm}
\usepackage{amscd}

\newcommand{\Q}{{\mathbb Q}}
\newcommand{\R}{{\mathbb R}}

\newcommand{\C}{{\mathbb C}}

\newcommand{\A}{{\mathbb A}}
\newcommand{\fg}{{\mathfrak g}}

\DeclareMathOperator*{\Injlim}{\varinjlim}

\newtheorem{theorem}{Theorem}
\newtheorem{corollary}{Corollary}

\theoremstyle{definition}

\numberwithin{conj}{section}
\newtheorem*{example}{Example}

\theoremstyle{remark}

\begin{document}

\title[Betti numbers of hyperbolic manifolds]{On the non-vanishing  
of the first Betti\\ number of hyperbolic three manifolds}
\date{}
\author{C.~S.~Rajan}

\address{Tata Institute of Fundamental 
Research, Homi Bhabha Road, Bombay - 400 005, INDIA.}
\email{rajan@math.tifr.res.in}

\subjclass{Primary 11F75; Secondary 22E40, 57M50}

\begin{abstract}
We show the
non-vanishing of  cohomology groups of sufficiently small congruence
lattices in $SL(1,D)$, where $D$ is a quaternion division
algebras defined over a number field $E$ contained
inside a solvable extension of a totally real number field. 
  As a corollary, we obtain
new examples of compact,  arithmetic,  hyperbolic three manifolds,
with non-torsion first homology group,  confirming a conjecture of
Thurston.   The proof uses the characterisation of the image of
solvable base change by the author, 
and the construction of cusp forms with non-zero
cusp cohomology by Labesse and Schwermer.

\end{abstract}

\maketitle

\section{Introduction}

Let $D$ be a quaternion division algebra over a number field  $E$.
Let $G$ denote the connected, semisimple algebraic group $SL_1(D)$
over $E$. Denote by   $G_{\infty}(E)$ the real Lie group $G(E\otimes
\R)$ and  fix a  maximal compact subgroup  $K_{\infty}$ of
$G_{\infty}$. Let $s_1$ (resp. $2r_2$) be the number of real
(resp. complex) places of $E$ at which $D$ splits,  and let
$s=s_1+r_2$.  The quotient space $M:=G_{\infty}/K_{\infty}$ with the
natural $G_{\infty}(E)$-invariant metric, is a symmetric space
isomorphic to ${\mathcal H}_2^{s_1}\times {\mathcal H}_3^{r_2}$, where
for a natural number $n$, ${\mathcal H}_n$ denotes the simply
connected  hyperbolic space  of dimension $n$.
  
Let $\A$ (resp. $\A_f$) denote the ring of adeles (resp. finite adeles) of
$\Q$.  Let $K$ be a compact, open subgroup of $G(\A_f\otimes E)$, and
denote by $\Gamma_K$ the  corresponding congruence arithmetic lattice
in $G_{\infty}(E)$ defined by the projection to $G_{\infty}(E)$ of the
group $G(E)\cap G_{\infty}(E)K$. For sufficiently small congruence
subgroups $K$, $\Gamma_K$ is a torsion-free lattice and
$\Gamma_K\backslash M$ is a (compact) Riemannian manifold.

In this note,  we prove
\begin{theorem} \label{constantcoeff}
With the above notation, assume further that  $E$ is a finite
extension of a totally real number field $F$ contained inside a
solvable extension $L$ of $F$.  For sufficiently small congruence
subgroups $\Gamma\subset G_{\infty}(E)$, the cohomology groups
\[H^s(\Gamma\backslash M, \C)\]
are non-zero.
\end{theorem}

The theorem was proved by Labesse and Schwermer \cite[Corollary
6.3]{LS}, in the case when there exists a tower of field extensions
\[ E=F_l\supset F_{l-1}\supset\cdots\supset F_0=F,\]
such that $F_{i+1}/F_i$ is either a cyclic extension of prime degree
or a non-normal cubic extension. The proof rests on the following two
observations: one, that the  base change (constructed by Langlands
\cite{L}) of the discrete series representations  from $SL(2,\R)$ to
$SL(2,\C)$ are  cohomologically
non-trivial representations of $SL(2,\C)$, and thus the base change of
cohomologically non-trivial cusp forms (for $SL(2)$) are
cohomologically non-trivial. 
 Secondly, when working with the group $SL_2$ over a
totally real field, the required cohomologically non-trivial
automorphic representations have discrete series as their archimedean
components, and an argument using pseudo-coefficients shows the
existence of such cusp forms. 

The new ingredient that goes into extending the theorem of Labesse and
Schwermer to all relatively solvable extensions of a totally real
field, is  the criterion of base change descent
for an invariant cuspidal automorphic representation with respect to a
solvable group of automorphisms of the field proved by the author in \cite{R}. 

If 
either $D$ is unramified at all finite places of
$E$ or if $F$ is taken to be the field of rationals and 
 the Galois closure $L$ of $E$ over $\Q$ is of odd
order over $\Q$, then Theorem \ref{constantcoeff}
was proved by Clozel \cite{C}. Clozel's proof uses the construction of
algebraic Hecke characters due to Weil, and the automorphic induction
of suitable such characters produces the desired cusp forms with
non-zero cohomology. Clozel's method proves the non-vanishing of the
cuspidal cohomology of the split groups $SL_2$ over any number
field $E$.  But in order to produce cohomological forms on inner forms of
$SL_2$ using the Jacquet-Langlands theorem, it is required that the
constructed cusp form on $SL_2(\A\otimes E)$ have discrete series
components at the places of $E$ where $D$ ramifies, and thus the
choice of $D$ has to be suitably restricted.

A particular case of interest is the following: 
\begin{corollary} \label{thurston}
With notation as in Theorem \ref{constantcoeff},   assume further that 
$E$ has exactly one pair of conjugate complex places, and the
quaternion division algebra $D$ is
ramified at all the real places of $E$. 
 For sufficiently small congruence subgroups $\Gamma\subset
G_{\infty}(E)$, the first betti number of the compact, hyperbolic
three manifold $\Gamma\backslash M$
is  non-zero.
\end{corollary}  
 A folklore conjecture (attributed to Thurston) is that the first
betti number of a compact, hyperbolic three manifold becomes positive
upon going to some finite cover. The first examples of compact,
hyperbolic arithmetic three manifolds $M_{\Gamma}$  with non-vanishing
rational first homology group are due to Millson \cite{M}. Using
geometric arguments,  Millson  showed the non-vanishing of the
first betti number for sufficiently small congruence subgroups, where
the arithmetic structure arises from rank $4$ quadratic forms over a
totally real number field $F$, and of signature $(3,1)$ at one
archimedean place and anisotropic at all other real places. As
mentioned above, other non-vanishing results confirming Thurston's
conjecture have been proved by Labesse-Schwermer
\cite{LS} and Clozel \cite{C}.

\begin{example} Let $L$ be a  non-real Galois extension of $\Q$ with Galois
group isomorphic to $S_4$. Consider a subgroup $H$ isomorphic to $S_3$
and containing a complex conjugation $\sigma$ corresponding to some
archimedean place. Then the invariant field $E=L^H$  is a
degree $4$ extension of $\Q$ having the following properties:
\begin{itemize}
\item $E$ does not contain any quadratic extension of $\Q$.
\item $E$ has at least one real place.
\item The subgroup containing the conjugates of $\sigma$ with respect
to $S_4$ will not be equal to $S_3$. Hence, $E$ is not totally real.
\end{itemize}
This gives us an example of a field $E$ satisfying the hypothesis of
the corollary. The non-vanishing result for the corresponding  hyperbolic three
manifolds given by Corollary \ref{thurston} are not in general 
covered by the non-vanishing results proved in 
\cite{C}, \cite{LS} or in \cite{M}.

\end{example}

\section{General coefficients} 

Theorem \ref{constantcoeff} can be generalized for suitable non-trivial
coefficient systems also. 
Let $F$ and $E$ be as in the hypothesis of the theorem. 
 Given a finite dimensional complex representation
$V$ of $SL_2(\R\otimes F)$, we now define the base change
representation $\Psi(V)$ of the group $G_{\infty}(E)$ \cite{LS}. We define
it first when $V$ is irreducible and extend it additively. If $V$ is
irreducible, then $V$ can be written as, 
\[ V\simeq \otimes_{v\in P_{\infty}(F)}V_{v}, \]
where $P_{\infty}(F)$ is the collection of the archimedean places of
$F$, and  the component $V_{v}$ of $V$ at the place $v$ is an 
 irreducible representation of
$SL_2(F_v)\simeq SL_2(\R)$, say  of dimension $k(v)$. 

Let $V_k$
(resp. $\bar{V}_k$) denote the irreducible, holomorphic
(resp. anti-holomorphic) representation of $SL_2(\C)$ of dimension
$k$. Restricted to $SU(2)$ they give raise to isomorphic
representations, which we continue to denote by $V_k$.  
Define the representation $W_k$ of $SL_2(\C)$ by $W_k=V_k\otimes
\bar{V}_k$. 

Suppose $D$ is a quaternion algebra over $E$.  We define the base
change coefficients $\Psi(V)$ of $G_{\infty}(E)$, as a tensor product
of the representations $\Psi(V)_w$ of the component groups $G(E_w)$,
as  $w$ runs over the collection of archimedean places of $E$. 
Suppose $w$ lies over a place $v$ of $F$. Define, 
\[ \Psi(V)_w\simeq \begin{cases} V_{k(v)} &  \text{if $w$ is real},\\
W_{k(v)} & \text{if $v$ is complex}. 
\end{cases}
\]

 Restricting the representation $\Psi(V)$ to a torsion-free  lattice
$\Gamma$ gives raise to  a 
well defined local system ${\mathcal L}_{\Psi(V)}$  on the
manifold $\Gamma\backslash M$.  The extension of Theorem
\ref{constantcoeff} to non-trivial coefficients  is the following:
\begin{theorem} \label{twistcoeff}
Let $F$ be a totally real number field, and  $L$ be a solvable finite
extension of $F$. Let $E$ be  a finite extension of $F$ contained in
$L$, and $D$ be a quaternion division algebra over $E$. Let $V$ be a
finite dimensional complex representation of $SL_2(\R\otimes
F)$. Then,
\[H^s(\Gamma\backslash M, {\mathcal L}_{\Psi(V)})\neq 0.\]
\end{theorem}

\section{Proof}
In order to prove Theorem \ref{twistcoeff}, it is more convenient to
work with  the cohomology groups $ H^*(G, E;V)$, defined  as a  direct
limit indexed by the compact open subgroups $K\subset G(\A_f\otimes
E)$:
\[  H^*(G, E;V)=\Injlim_{K}H^*(\Gamma_K, \Psi(V))\simeq
\Injlim_{K}H^*(\Gamma_K\backslash M, {\mathcal L}_V) . \] 
These cohomology groups
can be reinterpreted in terms of the automorphic spectrum of
$G(\A\otimes E)$.  Let $\rho$ denote the representation of
$G(\A\otimes E)$ acting by right translations on the space
$L^2(G(E)\backslash G(\A\otimes E))$ consisting of square integrable
functions  on $G(E)\backslash G(\A\otimes E)$. This decomposes as a
 direct sum of irreducible admissible representations $\pi$ of
$G(\A\otimes F)$ with finite multiplicity $m(\pi)$:
\[ \rho=\oplus_{\pi}m(\pi)\pi,\]
With respect to the decomposition $G(\A\otimes
E)=G_{\infty}(E)G(\A_f\otimes E)$, write $\pi=\pi_{\infty}\pi_f$,
where $\pi_{\infty}$ (resp. $\pi_f$) is a representation of
$G_{\infty}(E)$ (resp. $G(\A_f\otimes E)$).  The cohomology groups
$H^*(G, E;V)$ can also be expressed in terms of the relative Lie
algebra cohomology (see \cite{BW}) of the  automorphic spectrum as,
\begin{equation}\label{coh-autspectrum}
 H^*(G, E;V)\simeq \oplus_{\pi}m(\pi) H^*(\fg, K_{\infty},
\pi_{\infty}\otimes V)\otimes \pi_f,
\end{equation}
where $\fg$ is the Lie algebra of
$G_{\infty}(E)$. Hence in order to prove the theorem, it is enough to
construct an irreducible representation $\pi$ of $G(\A_E)$ with
$m(\pi)$ positive and such that $H^s(\fg, K_{\infty},
\pi_{\infty}\otimes V)$ is non-zero. 

We can assume that $V$ is irreducible of the form $ V\simeq
\otimes_{v\in P_{\infty}(F)}V_{k(v)}$. Let $D_k^+$ (resp. $D_k^-$) be
the holomorphic (resp. antiholomorphic) discrete series of $SL_2(\R)$
of weight $k+1$. We have, 
\begin{equation}\label{coh-ds}
H^q(\mathfrak{sl}_2(\R), SO(2), D_k^{\pm}\otimes V_k)=\begin{cases} 
\C & \text{if $q=1$},\\
0 &\text{otherwise},
\end{cases}
\end{equation}
where $\mathfrak{sl}_2(\R)$ and $\mathfrak{sl}_2(\C)$ denotes
respectively the Lie algebras of $SL_2(\R)$ and $SL_2(\C)$. 

Let $S$ be a finite set of finite places of $F$, 
containing all the finite places $v$ of $F$ dividing a 
finite place of $E$ at which $D$ ramifies.  By
\cite[Proposition 2.5]{LS}, there exists an irreducible, admissible
representation of $SL_2(\A\otimes F)$ satisfying the following
properties:
\begin{itemize} 
\item The multiplicity   $m_0(\pi)$ of $\pi$ occuring in the cuspidal
spectrum  $ L^2_0(SL_2(F)\backslash SL_2(\A\otimes F))$ consisting of
square integrable cuspidal functions on $SL_2(F)\backslash
SL_2(\A\otimes F)$
 is nonzero.  Further
$\pi$ is stable in the sense of \cite{LL}.
\item The local component $\pi_v$ of $\pi$ at an archimedean place $v$
of $F$ is a discrete series representation, with $\pi_{v}\in
\{D_{k(v)}^+,D_{k(v)}^-\}$. 
\item For any $v\in S$, the local component $\pi_v$ of $\pi$ is
  isomorphic to the Steinberg representation of $SL_2(F_v)$.
\end{itemize}
Let $\Pi$ be a cuspidal, automorphic representation of $GL_2(\A\otimes
F)$, such that $\pi$ occurs in the restriction of $\Pi$ to
$SL_2(\A\otimes F)$. Let $\Pi_{L}$ be the base change of $\Pi$ to
$GL_2(\A\otimes L)$ defined by Langlands in \cite{L}. Since $\pi$ is
stable, i.e., $\Pi$ is not automorphically induced from a character of
a quadratic extension of $F$, $\Pi_{L}$ is a cuspidal automorphic
representation of $GL_2(\A\otimes L)$.

Let $H$ be the Galois group of $L$ over $E$. Since $\Pi_{L}$ is
$H$-invariant, by the descent criterion proved in \cite{R}, there
exists an idele class character $\chi$ of $L$, such that the
representation $\Pi_{L}\otimes \chi$ is the base change from $E$ to
$L$ of a cuspidal representation $\Pi_E$ of $GL_2(\A\otimes E)$. Let
$\pi_E$ be a constituent of the restriction of $\Pi_E$ to
$SL_2(\A\otimes E)$, and occuring in the automorphic spectrum
$G(\A_E)$ with non-zero multiplicity $m(\pi_E)$.

 Base change makes sense at the level of $L$-packets (see \cite{LS}), 
and let $\pi_{k,
\C}$ denote the representation of $SL_2(\C)$ obtained as base change
of the $L$-packet $\{D_k^+, D_k^-\}$ ($L$-packets for complex groups
consist of only one element). It is known that (see \cite{LS}),
\begin{equation}\label{coh-sl2c}
H^1({\mathfrak{sl}}_2(\C), SU(2), \pi_{k,\C}\otimes W_k)\neq 0.
\end{equation}

 Let $w$ be an archimedean place of $E$ lying over a real place $v$ of
$F$.  Now twisting by a character does not alter the restriction of an
automorphic representation of $GL_2$ to $SL_2$. Hence if $w$ is a real
place of $E$, then the local component $\pi_{E,w}$ of $\pi_E$ at $w$
belongs to $\{D_{k(v)}^+,D_{k(v)}^-\}$, and if $w$ is a complex place
of $E$, then $\pi_{E,w}$ is isomorphic to $\pi_{k(v), \C}$. 

   The local components of the base change
to $E$ of $\pi$ continues to be the Steinberg representation of
$SL_2(E_w)$, at the places of $E$ where
$D$ ramifies. By the theorem of Jacquet-Langlands (\cite{JL},
\cite{LS}) applied
to $L$-packets of $SL_2$ and it's inner forms, we get an
automorphic representation $JL(\pi_E)$ of $G$ over $E$. At a place
$w$ where $D$ is ramified, the local component $JL(\pi_E)_w$ is
isomorphic to the restriction of the representation $V_{k(v)}$ to
$SU(2)$, where $v$ is a place of $F$ dividing $w$. In particular,  the zeroth
relative Lie cohomology group  
\begin{equation}\label{coh-su2}
H^0(\mathfrak{su}_2, SU_2, V_k\otimes V_k)=(V_k\otimes
V_k)^{SU(2)}\neq 0. 
\end{equation}
 At a  place $w$ of $E$ where $D$ splits,
$JL(\pi_E)_w\simeq \pi_{E,w}$, and hence the first  relative Lie algebra
cohomology with coefficients in the component of $\Psi(V)$ at $w$ is
non-zero.   It follows from equations (\ref{coh-ds}),
(\ref{coh-sl2c}), (\ref{coh-su2}) and by  the Kunneth
formula for the relative Lie algebra cohomology that
\[ H^{s}({\mathfrak g}, K_{\infty},
JL(\pi_E)_{\infty}\otimes \Psi(V))\neq 0.\]
By Equation (\ref{coh-autspectrum}),  this  proves Theorem \ref{twistcoeff}.


\begin{thebibliography}{JPSH}

\bibitem[BW]{BW} Borel, A. and Wallach, N. {\em Continuous cohomology,
discrete subgroups and representations of reductive groups},
Ann. Math. Stud. {\bf 94}, Princeton Univ. Press, 1980.

\bibitem[C]{C} Clozel, L.  {\em On the cuspidal cohomology of
arithmetic subgroups of ${\rm SL}(2n)$ and the first Betti number of
arithmetic $3$-manifolds},   Duke Math. J. {\bf 55} (1987), no. 2,
475--486.



\bibitem[JL]{JL} Jacquet, H. and Langlands, R. {\em Automorphic forms
on $GL(2)$}, Lect. Notes in Math. {\bf 114}, Berlin, Springer 1970.

\bibitem[LL]{LL} Labesse, J.-P. and Langlands, R. {\em
$L$-indistinguishability for $SL(2)$}, Can. J. Math. {\bf 31} (1979)
726-785.

\bibitem[LS]{LS} Labesse, J.-P. and Schwermer, J. {\em On liftings and
cusp cohomology of arithmetic groups}, Invent. math. {\bf 83} (1986)
383-401.


\bibitem[L]{L} Langlands, R. {\em Base change for GL(2)}, Annals of
Maths Studies {\bf 96} (1980), Princeton Univ. Press.

\bibitem[M]{M} Millson, J. J. {\em On the first Betti number of a
constant negatively curved manifold}, Ann. Math. {\bf 104} (1976)
235-247.

\bibitem[R]{R} Rajan, C. S.  {\em On the image and fibres of solvable
base change}, Math. Res. Letters, {\bf 9} (2002) no. 4,  499-508.




\end{thebibliography}
\end{document}